\documentclass[letterpaper, 10 pt, conference]{ieeeconf}  

\IEEEoverridecommandlockouts                              

\usepackage[utf8]{inputenc}
\usepackage{amssymb}
\usepackage{amsmath}
\usepackage{bm}
\usepackage{hyperref}
\usepackage{tikz}
\usepackage{tkz-euclide}
\usepackage{chngcntr}
\usepackage{verbatim}
\usepackage{mathtools}
\usepackage{esvect}
\usepackage{subfiles}
\usepackage{cite}
\usepackage{color}
\usepackage{xspace}
\usepackage{xparse}
\usepackage[noend]{algpseudocode}
\usepackage{subfigure}
\RequirePackage{algorithm}
\input{mysymbol.sty}

\usepackage{pgfplotstable}
\usepackage{pgfplots}
\pgfplotsset{
	table/search path={plot_figures},
}
\usepackage{tikz}
\usepackage{xr}
\usetikzlibrary{arrows,%
	petri,%
	topaths}%
\usetikzlibrary{graphs,graphs.standard}
\usepackage{tkz-berge}
\usepackage{booktabs}

\newtheorem{theorem}{Theorem}
\newtheorem{lemma}[theorem]{Lemma}

\newtheorem{remark}[theorem]{Remark}
\newtheorem{assumption}{Assumption}
\newtheorem{definition}{Definition}

\newcommand{\Oc}{\mathcal{O}}
\renewcommand{\reals}{\mathbb{R}}

\newcommand{\nag}{\textsc{Nag}\xspace}

\newcommand{\inner}[1]{\langle #1 \rangle}

\title{\LARGE \bf
Achieving Acceleration in Distributed Optimization\\ via Direct Discretization of the Heavy-Ball ODE
}

\author{Jingzhao Zhang, C\'esar A. Uribe, Aryan Mokhtari, and Ali Jadbabaie 
\thanks{The authors are with the Laboratory for Information and Decision Systems (LIDS), and the Institute for Data, Systems, and Society (IDSS),
        Massachusetts Institute of Technology, 77 Massachusetts Ave, Cambridge, MA 02139
        {\tt\small \{jzhzhang,cauribe,aryanm,jadbabai\}@mit.edu}}%
\thanks{*This research was supported in part by DARPA Lagrange and a Vannevar Bush Fellowship.}
}

\begin{document}

\maketitle

\begin{abstract}

\textcolor{black}{We develop a distributed algorithm for convex Empirical Risk Minimization,the problem of  minimizing large but finite sum of convex functions over networks. The proposed algorithm is derived from directly  discretizing the second-order heavy-ball differential equation and results in an accelerated convergence rate, i.e, faster than distributed gradient descent-based methods for strongly convex objectives that may not be smooth. Notably, we achieve acceleration without resorting to the well-known Nesterov's momentum approach.  We provide numerical experiments and contrast the proposed method with recently proposed optimal distributed optimization algorithms.}

\end{abstract}
\section{Introduction}
Acceleration in first-order optimization algorithms has recently become an intense focus of attention in machine learning, optimization, and related fields. From its original conception in the seminal works of Polyak~\cite{polyak1964some}, Nesterov~\cite{nes83}, and Nemirovskii~\cite{Nemirovskii1983}, the phenomena of acceleration has become central to the theory of convex optimization and its algorithms. This is mainly due to the fact that accelerated algorithms attain the theoretical oracle lower bounds for particular classes of convex problems~\cite{nes13}. Moreover, these lower bounds have been extended to the problem of distributed optimization over networks~\cite{sca17,scaman2018optimal}, where large quantities of data, as well as privacy constraints or limited access to complete information,  hinders the applicability of traditional centralized approaches~\cite{ned15}. Nevertheless, the fundamental understanding of acceleration remains a challenging problem. 

Recently, there have been many attempts to understand the acceleration phenomenon from a control theoretic viewpoint~\cite{all14,lessard2016analysis,fazlyab2018analysis}. In particular, many researchers have analyzed accelerated algorithms using a continuous-time interpretation suggesting that accelerated gradient method follows the trajectory of a second-order ordinary differential equation (ODE)~\cite{su-differential,wibisono-variational, diakonikolas2017accelerated,krichene2015accelerated}. Authors in \cite{zhang2018direct} further showed that direct discretization of a second-order ODE, which is also known as heavy ball ODE, generates accelerated first-order methods. In this paper, we aim to extend the discretization technique proposed in \cite{zhang2018direct} to the distributed setting. In particular, we design a novel accelerated first-order algorithm for solving a distributed optimization problem by  direct discretization of the heavy-ball ODE corresponding to the dual function of the problem.  

Our main result stipulates that the proposed algorithm achieves a convergence rate that is provably-faster than the one for gradient descent in terms of the number of communication rounds. Particularly, we show that the distance to a consensus solution among agents decreases at a sublinear rate of  $\mathcal{O}(N^{\frac{-2s}{s+1}})$ where $s$ is the order of the integrator used. This shows that by increasing the order of integrator the convergence rate of the proposed method approaches the optimal rate of $\mathcal{O}(N^{-2})$~\cite{uribe2018dual}. However, the experimental results provided in this work show that $s=4$ is sufficient to achieve a performance comparable with the optimal Nesteorv's accelerated gradient method. Moreover, the distance to the primal optimal solution is shown to decrease at a rate of $\mathcal{O}(N^{\frac{-2s}{s+1}})$ and the primal objective function suboptimality converges to zero at a rate of $\mathcal{O}(N^{\frac{-s}{s+1}})$.  

This paper is organized as follows. Section~\ref{sec:problem} recalls the problem of distributed optimization over networks and presents its formulation in terms of a set of equality constraints related to the Laplacian of the network as in~\cite{tutunov2016distributed}. Section~\ref{sec:dual} describes the dual formulation of the distributed optimization problem and some basic properties of the dual problem. Section~\ref{sec:algorithm} introduces the proposed algorithm and its derivation based on the discretization of the heavy-ball ODE corresponding to the dual formulation of the distributed optimization problem. Section~\ref{sec:convergence} presents our main analysis on the convergence properties of the proposed algorithm strongly convex functions that may not be smooth. Section~\ref{sec:numerical} shows numerical results of the proposed algorithm for two distributed optimization problems and compares the performance of our proposed method with its optimal counterparts. Finally in Section~\ref{sec:conclusions}, we summarize and discuss potential future work.

\vspace{2mm}
\noindent\textbf{Related work.}
Distributed consensus optimization has been studied heavily over the last decade. In particular, there have been several works which achieve a linear convergence rate for the setting that the local functions are strongly convex and smooth \cite{shi2014linear,makhdoumi2017convergence,shi2015extra,mokhtari64dqm,qu2017accelerated,ned17r}. For the case of convex and smooth functions, sublinear rates of $\mathcal{O}(N^{-1})$ have been proven using gradient descent based methods \cite{nedic2009distributed,jakovetic2014fast} as well as ADMM-type algorithms \cite{makhdoumi2017convergence,wei20131}. Perhaps the most related papers to our work are \cite{lan17} and  \cite{uribe2018dual} which study convex but non-smooth functions. \cite{lan17} approached the problem by regularizing the dual function and applying Nesterov's accelerated gradient method. \cite{uribe2018dual} used a similar idea but attained better dependency on the communication graph topology via a change of variable. Our approach follows a different path by discretizing the heavy-ball ODE defined by the dual objective. This approach does not require using regularization or Nesterov's accelerated gradient descent. Table ~\ref{tab:comparisons} shows a comparison between the communication cost of our proposed method and the distributed approaches introduced in \cite{lan17} and  \cite{uribe2018dual}. It further provides available convergence rates for centralized approaches in terms of gradient computations.
It is worth mentioning that convergence rates we refer in Table~\ref{tab:comparisons} for the distributed approaches, e.g., \cite{lan17} and  \cite{uribe2018dual}, are in terms of communication rounds and not gradient computations. Particularly, the dual approach assumes the availability of exact solutions to an auxiliary problem. Later in the paper we will discuss explicit dependencies on the function parameters and graph spectral properties for these algorithms.

\begin{table} \label{tab:comparisons}
\setlength{\tabcolsep}{2pt}

\begin{tabular*}{0.48\textwidth}{c|c|ccc}

\toprule
 & Centralized & \multicolumn{3}{c}{Decentralized} \\
  & (Gradient Comp.) & \multicolumn{3}{c}{(Communication Rounds)} \\
\midrule
Approach & \cite{Nemirovskii1983} & \cite{lan17} &  \cite{uribe2018dual} &  This work\\ 
\midrule
$\|Lx\|$ & - & $\Oc{(N^{-2})}$ &  $\Oc{(N^{-2})}$  &  $\Oc{(N^{-2\frac{s}{s+1}})}$\\ 
\midrule
$\|x_N - x^*\|^2$ & - & $\Oc{(N^{-2})}^\dagger$ &  $\Oc{(N^{-2})}^\dagger$ &  $\Oc{(N^{-2\frac{s}{s+1}})}$\\ 
\midrule
$f(x_N) - f(x^*)$ & $\Oc{(N^{-1})}$ & $\Oc{(N^{-2})}$ &  $\Oc{(N^{-2})}$ &  $\Oc{(N^{-\frac{s}{s+1}})}$\\ 
\bottomrule
\multicolumn{5}{l}{$\dagger$ Computed by $\|x_N - x^*\|^2 \le (f(x_N) - f(x^*))/\mu$}
\end{tabular*}
\caption{Iteration complexity of centralized and decentralized approaches for strongly convex non-smooth problems. }
\end{table}

\vspace{2mm}
\noindent\textbf{Notation.} 
For a matrix $A\in \reals^{m\times n}$, we denote $ker(A) = \{x \in \reals^n | Ax=0\}$ and $ker(A)^\perp = \{x \in \reals^n | x^Tv = 0\ \forall\ v \in ker(A)\}$. For a symmetric matrix $L$, we let $\lambda_{\max}(L)$ be its largest eigenvalue and $\lambda_{\min}^+(L)$ be its smallest positive eigenvalue. When $L$ is positive semidefinite, we further define $\sqrt{L}$ to be the unique positive semidefinite matrix such that $\sqrt{L}\sqrt{L} = L$. We
use $\boldsymbol{1}_n \in \reals^n$ and $\boldsymbol{0}_n\in \reals^n$ to denote vectors with all entries equal to 1 and 0, respectively. We will work in the standard Euclidean norm and let $\inner{\cdot, \cdot}$ denote its inner product.

\section{Problem Formulation}\label{sec:problem}

In this section, we formally define the distributed optimization problem over networks. Consider a set of $n$ nodes that communicate over a static, connected and undirected graph $\mathcal{G}=(\mathcal{V},\mathcal{E})$ where $\mathcal{V}=\{1,\cdots,n\}$ and $\mathcal{E} \subseteq \mathcal{V} \times \mathcal{V}$ denote the set of nodes and edges, respectively. We assume each node $i$ has access to a local convex function $f_{i}: \mathbb{R}^p\to  \mathbb{R}$, and nodes in the network cooperate to minimize the global objective function $f: \mathbb{R}^p\to  \mathbb{R}$ taking values $f(x)=\sum_{i=1}^n f_{i}(x)$. In other words, nodes aim to solve the optimization problem

\vspace{-0.3cm}
\begin{equation} \label{eq:main}
     \operatorname*{minimize}_{x \in \mathbb{R}^p}\ f(x) \ =\ \operatorname*{minimize}_{x \in \mathbb{R}^p}\ \sum_{i=1}^{n} f_i(x),
\end{equation}
while they are allowed to exchange information only with their neighbors. Two nodes $i$ and $j$ are considered neighbors, and therefore can communicate, if $(j,i) \in \mathcal{E}$. In this work, we assume that the local objective functions $f_i$ are strongly convex and as a result the global cost $f$ is also strongly convex.  Further, the local and global objective functions could be nonsmooth. As the objective function of Problem~\eqref{eq:main} is strongly convex, it has a unique solution denoted by $x^*$.

To solve Problem~\eqref{eq:main} in a decentralized fashion, the first step is assigning a local decision variable $x_i$ to each node $i$. Nodes aim to minimize the global objective function with their local information, while they ensure that their local decision variables are equal to their neighbors'. We use this interpretation to solve the following optimization problem

\begin{align} \label{eq:main_2}
   & \operatorname*{minimize}_{x_1,\dots,x_n \in \mathbb{R}^p}\ \sum_{i=1}^{n} f_i(x_i) \ \ \text{s.t.} \  x_i=x_j,\  \text{\forall}\ i,j \in \mathcal{V},
\end{align}
which is equivalent to Problem~\eqref{eq:main}, in the sense that the elements of a solution set $\{x_1^*,x^*_2,\dots,x^*_n\}$ of Problem~\eqref{eq:main_2} are equal  to the optimal solution of Problem~~\eqref{eq:main} which is $x^*$, i.e.,  $x_1^*=x^*_2=\dots=x^*_n=x^*$. 


We can simplify the notation in Problem~\eqref{eq:main_2} by defining  $\bbx=[x_1;\dots;x_n]\in \reals^{np}$ as the concatenation of the local decision variables $x_i$. Further, define $F: \mathbb{R}^{np}\to  \mathbb{R}$ as the sum of all local objective functions $F(\bbx)=F(x_1;\dots;x_n)=\sum_{i=1}^n f_i(x_i) $, and define matrix $L\in \reals^{n\times n}$ as the Laplacian matrix of the graph $\mathcal{G}$. It can be easily verified (see \cite{tutunov2016distributed}) that the constraint $x_1=\dots=x_n$ is equivalent to $\bbL\bbx=\bb0$ where $\bbL=L\otimes I_p\in \reals^{np}\times \reals^{np}$ is the Kronecker product of the Laplacian matrix $L$ and the identity matrix $I_p$. By incorporating these definitions Problem~\eqref{eq:main_2} can be written as
\begin{equation} \label{eq:prob}
    \min_{\bbx\in \reals^{np}}F(\bbx) \qquad \text{subject to}\ \sqrt{\bbL} \bbx=0.
\end{equation}
Note that for the constraint of \eqref{eq:prob} we can use $\sqrt{\bbL} \bbx=0$ instead of ${\bbL} \bbx=0$ as the matrix $\bbL$ is positive semidefinite and the null space of ${\bbL}$ and $\sqrt{\bbL}$ are identical. Hence, we can solve the matrix-form Problem~\eqref{eq:prob} in lieu of \eqref{eq:main_2} and the original problem in \eqref{eq:main}. Throughout  the rest of the paper, we use the notation $\bbx^*=[x^*;\dots;x^*]\in \reals^{np}$ to refer to the optimal solution of Problem~\eqref{eq:prob}.

\section{Dual Domain Representation}\label{sec:dual}

As projection to the null space of the matrix $\sqrt{\bbL}$ in a distributed fashion is not possible, to solve Problem \eqref{eq:prob} in the primal domain, one can minimize a penalized version of it; however, this approach yields convergence to a neighborhood of the optimal solution with a radius proportional to the penalty parameter \cite{yuan2016convergence}. One approach to designing a method with exact convergence is to solve Problem~\eqref{eq:prob} in the dual domain. In this work, we aim to solve the dual problem by discretizing its corresponding second-order heavy-ball ODE in a decentralized fashion.

We define the dual problem, which unlike the primal problem is unconstrained,  as 
\begin{equation}\label{dual_problem}
    \min_{\bby\in \reals^{np}}\ \varphi(\bby),
\end{equation}
where the dual function $\varphi(\bby)$ is defined as 
\begin{align} \label{eq:dual-fun}
    \varphi(\bby) = \max_{\bbx \in \reals^{np}} \left\{  \inner{\bby, \sqrt{\bbL}\bbx} - F(\bbx)\right\}.
\end{align}
Note that the dual function is convex, and due to the strong duality property the duality gap is zero. Further, the gradient of the dual function is given by 
\begin{align*}
    \nabla \varphi(\bby) = \sqrt{\bbL}\ \bbx^*\!\left(\!{\sqrt{\bbL}}\bby\!\right),
\end{align*}
where $\bbx^*(\bbz) : =\max_{\bbx} \left\{ \langle \bbz,\bbx \rangle -F(\bbx) \right\}$. This definition implies that evaluating the dual function gradient $\nabla \varphi(\bby)$ requires solving a convex program; however, in many cases this sub-problem either has a closed-form solution or can be solved efficiently. Functions for which one has immediate access to an explicit (or efficiently computed) solution $\bbx^*(\bbz)$ are sometimes called \textit{admissible} or \textit{dual-friendly}~\cite{Raginsky2012}.

In the rest of the section, we prove some properties of the dual function $\varphi$. Before doing so, however, we first formally state the required conditions in the following assumptions.

\vspace{2mm}
\begin{assumption} \label{assump:network}
The underlying communication graph $\ccalG$ is static, undirected and connected. 
\end{assumption}

\vspace{2mm}
\begin{assumption}\label{assump:mu}
The local objective functions $f_i$ are $\mu$-strongly convex, i.e., for any $x,y\in \reals^p$ 
\begin{equation*}
    f_i(y) \geq f_i(x) +  \inner{\nabla f_i(x), y-x} +\frac{\mu}{2} \|y-x\|^2.
\end{equation*}
\end{assumption}
\vspace{2mm}

For the next two assumptions, we need the requirements to be satisfied only over a compact set $\mathcal{B}\subset \reals^p$. Later in the paper (see \eqref{eq:x-compact-set}) we will properly define the set $\mathcal{B}$ which is determined by $\bbx^*$, $\lambda_{\max}(\bbL)$, $\lambda_{\min}^+(\bbL)$, $\min_{L\bbx=0}F(\bbx) - \min_\bbx F(\bbx)$ and $\mu$.

\vspace{2mm}
\begin{assumption}\label{assump:M}
The local objective functions $f_i$ are $M$-Lipschitz over the convex compact set $\mathcal{B}$, i.e., for any points $x,y\in \mathcal{B}$ 
\begin{equation*}
    |f_i(y) - f_i(x)| \le  M\|y-x\|.
\end{equation*}
\end{assumption}

\vspace{2mm}
\begin{assumption}\label{assump:s} $F(\bbx)$ is order $s+2$ differentiable on the compact set $\mathcal{B}$. 
\end{assumption}
\vspace{2mm}

The connectivity condition in Assumption~\ref{assump:network} implies that the Laplacian matrix $L$ satisfies 
\begin{equation*}
L=L^T,\ L\boldsymbol{1}_n=\boldsymbol{0}_n,\ \rank(L)=n-1.   
\end{equation*}
Further, it is easy to check that if the local functions $f_i$ are $\mu$-strongly convex and $M$-Lipschitz then the aggregated objective function $F$ is also $\mu$-strongly convex and $M$-Lipschitz. It is worth mentioning, if each function $f_i$ has a specific strong convexity $\mu_i$ and Lipschitz continuity $M_i$ parameters, then the aggregated objective function $F$ is $\min_{_i}\mu_i$ strongly convex and $\max_{_i}M_i$ Lipschitz.


 In the following lemma, we use differentiability and strong convexity of the function $F$ to show that the dual function $\varphi$ is differentiable.

\vspace{2mm}
\begin{lemma} \label{lemma:differentiability}
The dual function $\varphi(\bby)$ is $n$-th order differentiable at $\bby$ if $F(x)$ is $\mu-$strongly convex and $n$-th order differentiable at $\bbx^*(\sqrt{\bbL}\bby)$. Moreover, the dual function $\varphi(\bby)$ is $(\lambda_{\max}(\bbL)/\mu)$-smooth.
\end{lemma}

\vspace{2mm}
\begin{proof}
Recall that $\bbx^*(\sqrt{\bbL}\bby)$ is given by
\begin{equation} \label{eq:xstar}
    \bbx^*(\sqrt{\bbL}\bby) = \argmin_{\bbx\in\reals^{np}} \left\{\inner{\sqrt{\bbL}\bby, \bbx} - F(\bbx) \right\}.
\end{equation}
By KKT conditions, we know that 
\begin{equation} \label{eq:kkt-dual}
    \nabla \varphi(\bby) = \sqrt{\bbL} \bbx^*(\sqrt{\bbL}\bby), \quad
    \nabla F(\bbx^*(\sqrt{\bbL}\bby)) = \sqrt{\bbL} \bby.
\end{equation}
Moreover, we have 
\begin{equation*}
    \nabla^2 \varphi(\bby) = \sqrt{\bbL} \nabla_\bby \bbx^*, \quad
    \nabla^2 F(\bbx^*) \nabla_\bby \bbx^* = \sqrt{\bbL}.
\end{equation*}
By strong convexity, Hessian is invertible everywhere and $[\nabla^2 F(\bbx^*)]^{-1} \preceq I/\mu $. Therefore, the dual function is second-order differentiable and
\begin{equation} \label{eq:dual-hessian}
    \nabla^2 \varphi(\bby) = \sqrt{\bbL} [\nabla^2 F(\bbx^*)]^{-1} \sqrt{\bbL}.
\end{equation}
The desired result follows by recursively applying the following identity. Given a matrix valued function $K(x)$,
\begin{equation*}
    \nabla(K(\bbx))^{-1} = K(\bbx)^{-1}\nabla K(\bbx) K(\bbx)^{-1}.
\end{equation*}
In order to show that the dual function $\varphi$ is smooth, note that the eigenvalues of the primal Hessian inverse $ [\nabla^2 F(\bbx^*)]^{-1}$ are uniformly bounded above by $1/\mu$ due to strong convexity. Further, the eigenvalues of the square root of the Laplacian $\sqrt{\bbL} $ are also upper bounded by $\sqrt{\lambda_{\max}(\bbL)}$. Therefore, based on the expression in \eqref{eq:dual-hessian}, the eigenvalues of the dual function Hessian $\nabla^2 \varphi(\bby)$ are bounded above by $\lambda_{\max}(\bbL)/\mu$, and, hence, the dual function is $(\lambda_{\max}(\bbL)/\mu)$-smooth; see also 
 Proposition~$12.60$ in~\cite{rockafellar2011variational}.
\end{proof}


\section{Algorithm}\label{sec:algorithm}
In this section, we first review Runge-Kutta (RK) integrators. Then, we state our dynamical system of interest and introduce a distributed accelerated method by discretizing the dynamical system using Runge-Kutta integrators. 

\subsection{Runge-Kutta integrators}
Here, we briefly recap explicit Runge-Kutta (RK) integrators used in our work. For a more in-depth discussion, please see the textbook~\cite{hairer-textbook}. Consider a dynamical system $\dot{\bbzeta} = G(\bbzeta)$ and let the current point be $\bbzeta_0$ and the step size be $h$. An \emph{explicit $S$ stage Runge-Kutta method} generates the next step via the following update:
\begin{equation}
  \begin{split}\label{eq:runge-kutta}    
    g_i &= \bbzeta_0 + h \textstyle \sum_{j=1}^{i-1} a_{ij} G(g_j), \\
    \Phi_h(\bbzeta_0) &= \bbzeta_0 + h \textstyle \sum_{i=1}^{S} b_i G(g_i),
  \end{split}
\end{equation}
where $a_{ij}$ and $b_i$ are suitable coefficients defined by the integrator;  $\Phi_h(\bbzeta_0)$ is the estimation of the state after time step $h$, while $g_i$ (for $i=1,\ldots,S$) are a few neighboring points where the gradient information $G(g_i)$ is evaluated. 

By combining the gradients at several evaluation points, the integrator can achieve higher precision by matching up Taylor expansion coefficients. In the following definition we formally define the order of an integrator. 
\vspace{2mm}
\begin{definition}
Let $\Psi_h(\bbzeta_0)$ be the true solution to the ODE with initial condition $\bbzeta_0$ and $\Phi_h(\bbzeta_0)$ be the estimation of the state after time step $h$; we say that an integrator $\Phi_h(\bbzeta_0)$ has order $s$ if its \emph{discretization error} shrinks as
\begin{equation*}
\label{eq:disc-error}
\|\Phi_h(\bbzeta_0) - \Psi_h(\bbzeta_0)\| = O(h^{s+1}), \quad\text{as}\ h\to 0.
\end{equation*}
\end{definition}
\vspace{2mm}

In general, RK methods offer a powerful class of numerical integrators, encompassing several basic schemes. The explicit Euler's method defined by $\Phi_h(\bbzeta_0) = \bbzeta_0 + hG(\bbzeta_0)$ is an explicit RK method of order 1, while the midpoint method $\Phi_h(\bbzeta_0) = \bbzeta_0 + hG(\bbzeta_0+\tfrac{h}{2}G(\bbzeta_0))$ is of order 2. Depending on the order of the RK integrators the number of required stages, i.e., the number of gradient evaluations, varies. For instance, a Runge-Kutta integrator of degree $s$ where $1\leq s\leq 4$ requires $S=s$ stages, while a Runge-Kutta integrator of order $s=9$ requires $S=16$ stages , i.e., $16$ gradient evaluations per iteration.

\subsection{Optimization methods as dynamical systems}
We start with Nesterov's accelerated gradient  (\nag) method \cite{nes13} for convex smooth problems. To solve the dual problem using \nag we need to follow the updates
\begin{equation}
\label{eq:nesterov}
\begin{split}
\bby_k &= \bbz_{k-1} - h\nabla \varphi(\bbz_{k-1}),\\ 
\quad \bbz_k &= \bby_k + \tfrac{k-1}{k+2}(\bby_k - \bby_{k-1}),
\end{split}
\end{equation}
where $\bbz\in \reals^{np}$ is an auxiliary variable, $h$ is a positive stepsize, and $k$ is the iteration index. In \cite{su-differential}, the authors showed that the iteration in~\eqref{eq:nesterov} in the limit when $h\to 0$ is equivalent to the following second-order ODE
\begin{align}\label{eq:nag-ode}
\ddot{\bby}(t) + \tfrac{3}{t}\dot{\bby}(t) + \nabla \varphi(\bby(t)) = 0,  \qquad \text{where}\ \dot{\bby} = \tfrac{d \bby}{d t}.
\end{align}
This ODE is also known as heavy-ball ODE which relates to the heavy-ball method proposed by Polyak~\cite{polyak1964some}. 

It can be shown that in the continuous domain the objective function suboptimality gap  $\varphi(\bby(t))-\varphi(\bby^*)$ decreases at a rate of $\Oc({1}/{t^2})$ along the trajectories of the ODE (see \cite{su-differential, wibisono-variational}). The work in  \cite{zhang2018direct}, studied the reverse problem of discretizing \eqref{eq:nag-ode} to get stable optimization algorithms. In particular, it proposed a variation of the second-order ODE in (\ref{eq:nag-ode}) which can be written as the following dynamical system
\begin{align}\label{eq:dynamics}
\dot{\bbzeta} = 
G(\bbzeta)
= 
\begin{bmatrix}
-\frac{5}{t} \bbv - 4 \nabla \varphi(\bby)\\
\bbv\\
1
\end{bmatrix},\quad \bbzeta = [\bbv; \bby; t],
\end{align}
where the variable $\bbzeta\in \reals^{2np+1}$ is the concatenation of the decision variable $\bby$, its time derivative $\bbv=\tfrac{d \bby}{d t}$, and time $t$. It has been shown that direct discretization of the dynamical system in \eqref{eq:dynamics} with any explicit Runge-Kutta integrator leads to a stable algorithm \cite{zhang2018direct}. In particular, if the function $\varphi(\bby)$ is order $s+2$ differentiable with bounded high order derivatives (i.e., $\exists L$, s.t. $\forall \ p=2,3,...,s+2, \|\nabla^{(p)}\varphi(y)\| \le L $), discretizing the ODE with an order-$s$ Runge-Kutta integrator achieves a convergence rate of $\Oc{(N^{-2s/(s+1)})}$, which is faster than the $\Oc{(N^{-1})}$ convergence rate of gradient descent.

\begin{algorithm}[ht]
    Input ($f, M, \mu, s, a_{ij}, b_j$) \Comment{$s, a_{ij}, b_j$ are defined by the chosen Runge-Kutta integrator}
    \begin{algorithmic}[1]
        \State For each agent i
        \State Set the initial variable $\bbzeta_i = [\bb0; \bb0; 1]$
        \State Set step size $h = h_0/N^{\frac{s}{s+1}}$ \Comment{$h_0$ depends on $\mu, M, s$}
        \For {$k = 1, \dots , N$} 
        \For {$l = 1,\dots , S$}
        \State $\hat{g}_{l, k}^i = \xi_{k}^i + h \sum_{j=1}^{l-1} a_{lj} G(\hat{g}_{l, k}^i)$
        \State Denote components of $\hat{g}_{l, k}^i$ as $[\hat{v}_{l, k}^i; \hat{y}_{l, k}^i; t_{l, k}^i]$
        \State $x^*_{i, l}(\hat{y}_{l, k}^i) = \arg\max_{x_i} \inner{\hat{y}_{l, k}^i, x_i} - f_i(x_i)$
        \State Broadcast $x^*_{i, l}(\hat{y}_{l, k}^i)$ to neighbors
        \State $G(\hat{g}_{l, k}^i) = 
        \begin{bmatrix}
        -\frac{5}{\Tilde{t}_{l, k}^i} \hat{v}_{l, k}^i - 4 \sum_{j=1}^m L_{ij} x^*_{i, l}(\hat{y}_{l, k}^i)\\
        \hat{v}_{l, k}^i\\
        1
        \end{bmatrix}$
        \EndFor
        \State $ \xi_{k+1}^i = \xi_{k}^i + h \sum_{j=1}^{S}b_j G(\hat{g}_{j, k}^i)$
        \EndFor
    \end{algorithmic}
    \caption{Distributed optimization Based on direct discretization }\label{alg_1}
\end{algorithm}

\subsection{Distributed ODE discretization}
In this section, we propose a novel algorithm that solves Problem \eqref{eq:prob} in a \textit{decentralized} fashion by following the updates defined based on the direct discretization of the dynamical system in \eqref{eq:dynamics} using Runge-Kutta (RK) integrators. The sequence of iterates generated by RK discretization of the dynamical system \eqref{eq:prob} can be written as 
\begin{equation}\label{eq:runge-kutta-2}    
    \bbzeta_{k+1}= \bbzeta_k + h \textstyle \sum_{i=1}^{S} b_i G(\bbg_i),
\end{equation}
where $ \bbg_i = \bbzeta_k + h \textstyle \sum_{j=1}^{i-1} a_{ij} G(\bbg_j)$ and $G$ is given by
\begin{align}\label{eq:dynamics_2}
G(\bbzeta)
= 
\begin{bmatrix}
-\frac{5}{t} \bbv - 4 \sqrt{\bbL}\bbx^*\left({\sqrt{\bbL}}\bby\right)\\
\bbv\\
1
\end{bmatrix}.
\end{align}
Notice that the sparsity pattern of $\sqrt{\bbL}$ may be different from $\bbL$ and hence the operation $\sqrt{\bbL} \bbx$ cannot be executed over a network by exchanging information only with neighboring nodes. Therefore, we apply a change of variable $\bbxi := [\hat{\bbv}, \hat{\bby}, t] := [\sqrt{\bbL}\bbv, \sqrt{\bbL}\bby, t]$. Then the update step of Runge-Kutta integrator defined in \eqref{eq:runge-kutta-2} becomes
\begin{equation}
  \begin{split}\label{eq:runge-kutta_3}    
    {\bbg_{i,k}} &= \bbxi_k + h \textstyle \sum_{j=1}^{i-1} a_{ij} \hat{G}(\bbg_{j,k}), \\
    \bbxi_{k+1} &= \bbxi_k + h \textstyle \sum_{i=1}^{S} b_i \hat{G}(\bbg_{i,k}),
  \end{split}
\end{equation}
with the revised dynamical system $\hat{G}$ defined as 
\begin{align} \label{eq:runge-kutta_4}
\hat{G}(\bbxi)
= 
\begin{bmatrix}
-\frac{5}{t} \hat{\bbv} - 4 {\bbL}\bbx^*\left(\hat{\bby}\right)\\
\hat{\bbv} \\
1
\end{bmatrix}
\end{align}
for $\bbxi=[\hat{\bbv}, \hat{\bby}, t]$. Recall that, as defined in \eqref{eq:xstar}, the variable $\bbx^*(\hat{\bby})$ is given by $\bbx^*(\hat{\bby}) = \argmin_{\bbx\in\reals^{np}} \left\{\inner{\hat{\bby}, \bbx} - F(\bbx) \right\}$. The above dynamics can be evaluated in a distributed manner by exchanging information only with neighboring nodes as the graph Laplacian has the sparsity pattern of the graph $\mathcal{G}$. In particular, to perform the system of updates in \eqref{eq:runge-kutta_3} node $i$ can update its concatenated local decision variable $\xi_k^i=[\hat{v}_{k}^i;\hat{y}_k^i;t_k^i]\in \reals^{2p+1}$ at step $k$ based on the update
\begin{equation}\label{eq:local_update}
 \xi_{k+1}^i = \xi_{k}^i + h \sum_{j=1}^{S}b_j \hat{G}^i(\hat{g}_{j, k}^i),
\end{equation}
where the vectors $\hat{g}_{l, k}^i$ are defined as
\begin{equation*} 
\hat{g}_{l, k}^i = \xi_{k}^i + h \sum_{j=1}^{l-1} a_{lj} \hat{G}^i(\hat{g}_{l, k}^i)
\end{equation*}
and the operator $\hat{G}^i$ is given by
\begin{equation*} 
\hat{G}^i(\hat{g}_{l, k}^i) = 
        \begin{bmatrix}
        -\frac{5}{\Tilde{t}_{l, k}^i} \hat{v}_{l, k}^i - 4 \sum_{j=1}^m L_{ij} x^*_{i, l}(\hat{y}_{l, k}^i)\\
        \hat{v}_{l, k}^i\\
        1
        \end{bmatrix}.
\end{equation*}
The detailed steps of the proposed method are summarized in Algorithm~\ref{alg_1}. Note that the initial variables $y_0^i$ are set to $0_{p}$ so that $\bby_0 = \bb0_{np}$. This condition is needed to ensure that the sequence of variables $\bby_k$ are always orthogonal to the kernel space of the Laplacian matrix $\bbL$ as we show in Lemma~\ref{lemma:subspace}.

\begin{remark}
Each iteration of the proposed algorithm requires $S$ rounds of communications between neighboring nodes per iteration, as at each iteration each node $i$ has to evaluate $G(\hat{g}_{l, k}^i)$ for $S$ different points. 
\end{remark}

\section{Convergence analysis}\label{sec:convergence}

In this section, we state the theoretical convergence guarantees for our proposed algorithm. We further compare them against known optimal rates. To do this, we first prove the following auxiliary lemma.
\vspace{1mm}

\begin{lemma} \label{lemma:subspace}
If the initial dual variable is $\bby_0 = 0$, then for all $k > 0$, the dual variables $\bby_k$ are orthogonal to the kernel space of the Laplacian matrix $\bbL$, i.e., $\bby_k\in ker(\bbL)^{\perp}$. 
\end{lemma}

\vspace{1mm}
\begin{proof}
 To prove that $\bby_k$ is within the span of the square root of the Laplacian $\sqrt{\bbL}$ we first need to show that $\bbv_k$ satisfies this condition. According to the update of  $\bbv_k$ in \eqref{eq:runge-kutta-2} and the dynamical system in \eqref{eq:dynamics_2}, if we set the initial vector $\bbv_0=0$, then $\bbv_k$ is a linear combination of a set of vectors that can be written as $\sqrt{\bbL}(\bbx^*(\sqrt{\bbL}\bby))$ which are in the span of $\sqrt{\bbL}$. Then, based on the initial condition $\bby_0=0$ and the update of $\bby_k$ which only depends on $\bbv$, we obtain that $\bby_k$ is in the span of $\sqrt{\bbL}$. Therefore, $\bby_k\in \{\bbu| \bbu \perp ker(\bbL)\}$.
\end{proof}
\vspace{1mm}

Lemma~\ref{lemma:subspace} shows that the dual iterates $\bby_k$ always stay in the span of the Laplacian matrix $\bbL$. In the following theorem, We use this result to characterize the convergence guarantees of our proposed algorithm. 
\vspace{1mm}

\begin{theorem}\label{main_theorem}
Consider the proposed method outlined in Algorithm~\ref{alg_1}. Suppose that the conditions in Assumptions \ref{assump:network}- \ref{assump:s} are satisfied. Further, let $\bby_N$ be the dual iterate generated after running Algorithm~\ref{alg_1} for $N$ iterations using an \mbox{order-$s$} Runge-Kutta integrator with $S$-stages. Then, the primal variable $\bbx_N=\bbx^*(\sqrt{\bbL}\bby_N)$ corresponding to the iterate  $\bby_N$ satisfies the following inequalities:\\
(i) consensus distance
\begin{equation}\label{main_theorem_first_claim}
    \|\sqrt{\bbL}\bbx_N\|^2  \le \Oc{\left(\tfrac{\lambda_{\max}(\bbL)^3}{\mu^3}SN^{\frac{-2s}{s+1}}\right)},
\end{equation}
(ii) average distance to primal optimum
\begin{align}\label{main_theorem_second_claim}
    \frac{1}{n}\|\bbx_N - \bbx^*\|^2  \le \Oc{\left(\tfrac{\lambda_{\max}(\bbL)^2}{n \mu^3}SN^{\frac{-2s}{s+1}}\right)},
\end{align}
(iii) average aggregated objective suboptimality 
\begin{equation}\label{main_theorem_third_claim}
\frac{1}{n}[F(\bbx_N) \!-\! F(\bbx^*)] \le \Oc\!\left(\sqrt{\tfrac{S\lambda_{\max}(\bbL)^3}{n^2\mu^3\lambda_{\min}^+(\bbL)}}MN^{\frac{-s}{s+1}}\!\right)\!.
\end{equation}
\end{theorem}

\vspace{2mm}
\begin{proof}
 We prove the claims in four steps. In this proof, to simplify the notation, we denote $\bbx^*(\sqrt{\bbL}\bby_k)$ by $\bbx_k$. 

\emph{\textbf{Step 1: Show that $\{\bby_k\}_{k\ge 0}$ and $\{\bbx_k\}_{k\ge 0}$ stay in bounded sets.}} The boundedness of $\{\bbx_k\}_{k\ge 0}$ show that it suffices for Assumption 3 and 4 to hold on a bounded set, and The boundedness of $\{\bby_k\}_{k\ge 0}$ is needed for using the results in~\cite{zhang2018direct}.

Note that the difference $\bby_k-\bby^*$ can be written as $\bby_k-\bby^*= (\bby_k-\frac{t}{4}\bbv_k-\bby^*) +(\frac{t}{4}\bbv_k) $. Therefore, using the inueqality $\|a+b\|^2\leq 2\|a\|^2+2\|b\|^2$ we can write 
\begin{align} \label{eq:bound-y}
    &\|\bby_k - \bby^*\|^2 \le 2\|\bby_k - \frac{t}{4}\bbv_k - \bby^*\|^2 + \frac{t^2}{8}\|\bbv_k\|^2,
\end{align}
By Proposition 7 of  \cite{zhang2018direct}, after applying discretization for $k$ iterations ($k<N$), we have
\begin{align}\label{eq:proof_100}
    &\|\bby_k - \frac{t}{4}\bbv_k - \bby^*\|^2 + \frac{t^2}{16}\|\bbv_k\|^2,  \nonumber\\
    &\le \exp(1) (\| \bby^*\|^2 + \varphi(0)-\varphi(\bby^*)) + 1.
\end{align}
By combining \eqref{eq:bound-y} and \eqref{eq:proof_100} we obtain that 
\begin{align} \label{eq:proof_200}
    &\|\bby_k - \bby^*\|^2 \le  2\exp(1) (\| \bby^*\|^2 + \varphi(0)-\varphi(\bby^*)) + 2.
\end{align}
 According to the KKT condition of Problem~\eqref{eq:prob}, we can write $\nabla F(\bbx^*)+\sqrt{\bbL}\bby^*=0$. This implies  $\| \bby^*\|^2 \le \|\nabla F(\bbx^*)\|^2/\lambda_{\min}^+(\bbL)$. Using this inequality and the fact that $\varphi(0) = -\min_\bbx F(\bbx)$ and $\varphi(\bby^*) = -\min_{\bbL\bbx=0}F(\bbx)$ as well as the result in \eqref{eq:proof_200} we can derive the following bound 
\begin{align} \label{eq:proof_300}
    &\|\bby_k - \bby^*\|^2 \\
    &\le  2\exp(1) \left[\frac{\|\nabla F(\bbx^*)\|^2}{\lambda_{\min}^+(\bbL)} + \min_{L\bbx=0}F(\bbx) - \min_\bbx F(\bbx))\right] + 2.\nonumber
\end{align}
Hence, for all iterates $k$ the distance $\|\bby_k - \bby^*\|^2$ is bounded by a constant and the iterates  $\{\bby_k\}_{k\ge 0}$ lie in a compact bounded set defined as $\{\bby \in \reals^{np}| \|\bby_k - \bby_*\|^2 \le 2\mathcal{E} \}$, where
\begin{equation}\label{eq:proof_400}
    \mathcal{E} = \exp(1) \left[\frac{\|\nabla F(\bbx^*)\|^2}{\lambda_{\min}^+(\bbL)} + \min_{L\bbx=0}F(\bbx) - \min_\bbx F(\bbx)\right] + 1.
\end{equation}

Now, we proceed to show that the primal iterates stay within a bounded set. Using the strong convexity of the primal function $F$ we can write 
\begin{align}
    \|\bbx_k - \bbx^*\| \le \frac{\|\nabla F(\bbx_k) - \nabla F(\bbx^*)\|}{\mu} 
    \end{align}
    Now using the definition of the primal gradient in \eqref{eq:kkt-dual} we can replace $\nabla F(\bbx_k) - \nabla F(\bbx^*)$ by $\sqrt{\bbL} \bby_k-\sqrt{\bbL} \bby^*$ to obtain
    \begin{align*}
    \|\bbx_k - \bbx^*\| & \leq  \frac{\|\sqrt{\bbL} \bby_k-\sqrt{\bbL} \bby^*\|}{\mu} \le \frac{\sqrt{\lambda_{\max}(\bbL)}}{\mu}\|\bby_k - \bby^*\|,
\end{align*}
where the last inequality follows by the bound on the eigenvalues of the Laplacian. Therefore, based on the result in \eqref{eq:proof_300}, all primal vectors $ \{\bbx^*(\sqrt{L}\bby_k)\}_{k\geq0}=\{\bbx_k\}_{k\geq0} $ stay in a compact set $\mathcal{B}$ defined as
\begin{equation} \label{eq:x-compact-set}
\mathcal{B} : = \left\{\bbx\ \bigg{|}\ \|\bbx-\bbx^*\| \le \frac{\sqrt{2\mathcal{E}\lambda_{\max}(\bbL)}}{\mu}\right\},
\end{equation}
where $\mathcal{E}$ defined in \eqref{eq:proof_400} is a constant determined by $\lambda_{\min}^+(\bbL)$, $\min_{\bbL\bbx=0}F(\bbx) - \min_\bbx F(\bbx))$ and $\mu$. Note that this result shows that the iterates $\bbx_k$ stay in a bounded set $\mathcal{B}$ and therefore the conditions in Assumptions 3 and 4 should only hold over the compact convex set $\mathcal{B}$.

\emph{\textbf{Step 2: Apply Theorem 1 in \cite{zhang2018direct} to get nonasymptotic rate for the dual function.}} By Lemma \ref{lemma:differentiability}, we know that $\varphi$ is order $s+2$ differentiable. By continuity, we know that all its high order derivatives have bounded operator norm. In addition, we know that the dual iterates $\bby_k$ stay in a bounded set. 
Furthermore, we know $\varphi(\bby)$ has $\tfrac{\lambda_{\max}(\bbL)}{\mu}-$Lipschitz gradients. Hence, all the required conditions of Theorem~$1$ in \cite{zhang2018direct} are satisfied and we obtain that 
\begin{equation}\label{dual_error_bound}
   \varphi(\bby_N)-\varphi(\bby^*) = \Oc{\left(\tfrac{\lambda_{\max}(\bbL)^2}{\mu^2}SN^{\frac{-2s}{s+1}}\right)},
\end{equation}
where $S$ is the number of stages in the RK integrator.

\emph{\textbf{Step 3: Bound the distance to optimum of the primal objective.}} Recall the definition of the dual function in \eqref{eq:dual-fun}. As $\bbx_N= \bbx^*(\sqrt{\bbL}\bby_N)$ we know that
\begin{equation}\label{proof_210}
    \varphi(\bby_N) = \inner{\sqrt{\bbL} \bby_N, \bbx_N }-F(\bbx_N)
\end{equation}
\begin{equation}\label{proof_230}
    \varphi(\bby^*) = \inner{\sqrt{\bbL} \bby^*, \bbx^* }-F(\bbx^*)=-F(\bbx^*) ,
\end{equation}
where we used the fact that $\sqrt{\bbL}\bbx^* = 0$. Subtract the sides of \eqref{proof_230} from the ones in \eqref{proof_210} to obtain
\begin{align} \label{eq:primal-subop1}
    \varphi(\bby_N)\!-\!\varphi(\bby^*) & =\! \inner{\sqrt{\bbL} \bby_N, \bbx_N}\! -\! F(\bbx_N) + F(\bbx^*) \nonumber\\
    & =\! \inner{\sqrt{\bbL} \bby_N, \bbx_N \!-\! \bbx^*} \!-\! F(\bbx_N)\! +\! F(\bbx^*),
\end{align}
where in the second equality we again use $\sqrt{\bbL}\bbx^*=0$. Also, using $\mu$-strong convexity of the primal function $F$ and the expression in \eqref{eq:kkt-dual} we can write
\begin{align} \label{eq:primal-subop11}
    \|\bbx_N \!-\! \bbx^*\|^2 &\le \frac{2}{\mu}(\inner{\nabla F(\bbx_N), \bbx_N - \bbx^*} - F(\bbx_N) + F(\bbx^*)) \nonumber\\
    & = \frac{2}{\mu}(\inner{\sqrt{\bbL} \bby_N, \bbx_N\! -\! \bbx^*}\!-\! F(\bbx_N) \!+\! F(\bbx^*)).
\end{align}
Combining the results in \eqref{eq:primal-subop1} and \eqref{eq:primal-subop11} leads to
\begin{align}
    \|\bbx_N - \bbx^*\|^2 &\le \frac{2}{\mu}\left( \varphi(\bby_N)-\varphi(\bby^*)\right),
\end{align}
and hence by using the result in \eqref{dual_error_bound} we obtain that 
\begin{align}
    \|\bbx_N - \bbx^*\|^2 \le \Oc{\left(\tfrac{\lambda_{\max}(\bbL)^2}{\mu^3}SN^{\frac{-2s}{s+1}}\right)},
\end{align}
and the claim in \eqref{main_theorem_second_claim} follows. 

\emph{\textbf{Step 4: Bound the consensus distance and the suboptimality of the primal objective.}}
By Lemma~\ref{lemma:differentiability}, $\varphi(\bby)$ is $(\lambda_{\max}(\bbL)/\mu)$-smooth. By smoothness property (see \cite{nes13}),
\begin{equation}\label{eq:proof_bound_on_dual_norm}
    \|\nabla \varphi(\bby_N)\|^2 \le \frac{\lambda_{\max}(\bbL)}{\mu}(\varphi(\bby_N) - \varphi(\bby^*)).
\end{equation}
By using the definition of the dual gradient in \eqref{eq:kkt-dual} we can replace $ \|\nabla \varphi(\bby_N)\|^2$ in \eqref{eq:proof_bound_on_dual_norm} by $\|\sqrt{\bbL}\bbx_N\|^2$. Further, we can substitute $\varphi(\bby_N) - \varphi(\bby^*)$ in \eqref{eq:proof_bound_on_dual_norm} by its upper bound in \eqref{dual_error_bound} to obtain
\begin{align}\label{first_result}
    \|\sqrt{\bbL}\bbx_N\|^2 \leq  \Oc{\left(\tfrac{\lambda_{\max}(\bbL)^3}{\mu^3}SN^{\frac{-2s}{s+1}}\right)},
\end{align}
which yields the claim in \eqref{main_theorem_first_claim}.


 We then decompose $\bbx_N = \bbx_N^{\perp} + \bbx_N^\parallel$, where $\bbx_N^\parallel \in ker(\bbL)$ and $\bbx_N^\perp \perp ker(\bbL)$. Based on the result in \eqref{first_result} we can write
\begin{equation}
     \|\sqrt{\bbL}\bbx_N^\perp \|^2  = \|\sqrt{\bbL}\bbx_N\|^2 \leq \Oc{\left(\tfrac{\lambda_{\max}(\bbL)^3}{\mu^3}SN^{\frac{-2s}{s+1}}\right)}.
\end{equation}
By the fact that $\bbx_N^{\perp} \perp ker(\bbL)$, we have
\begin{equation} \label{eq:xperp}
    \|\bbx_N^\perp\|^2 \leq \Oc{\left(\tfrac{\lambda_{\max}(\bbL)^3}{\mu^3\lambda^{+}_{\min}(\bbL)} SN^{\frac{-2s}{s+1}}\right)}.
\end{equation}
Further, using the expression in \eqref{eq:xstar}, we can write
\begin{equation}\label{mid_step}
    \inner{\sqrt{\bbL}\bby_N, \bbx_N} - F(\bbx_N) \ge \inner{\sqrt{\bbL}\bby_N, \bbx^*}-F(\bbx^*)
\end{equation}
as $\bbx_N$ is the maximizer of the problem for $\bby_N$.
Rearrange the terms in \eqref{mid_step} and use the fact that $\sqrt{\bbL}\bbx^*=0$ to obtain
\begin{align}\label{close_to_final}
      F(\bbx_N) &\le F(\bbx^*) + \inner{\sqrt{\bbL}\bby_N, \bbx_N} \nonumber \\
      &\le F(\bbx^*) + |\inner{\sqrt{\bbL}\bby_N, \bbx_N}| \nonumber \\
      &= F(\bbx^*) + |\inner{\nabla F(\bbx_N), \bbx_N^\perp}|,
\end{align}
where the last inequality follows from ~\eqref{eq:kkt-dual} and $\sqrt{\bbL}\bbx_N^\parallel=0$. Using the Cauchy-Schwartz inequality and \eqref{close_to_final} we obtain
\begin{align} \label{eq:subopt-Cauchy}
      F(\bbx_N) -F(\bbx^*) \le   \|\nabla F(\bbx_N)\| \|\bbx_N^\perp\|,
\end{align}
Assumption~\ref{assump:M} implies that $\nabla \|F(\bbx)\|$ is bounded above by $M$ for any $\bbx$ in the set $\mathcal{B}$. As we showed in step 1 of the proof, the iterate $\bbx_N$ is within the set $\mathcal{B}$ and hence 
\begin{align}\label{upper_bound_on_primal_grad_norm}
    \|\nabla F(\bbx_N)\| \le M.
\end{align}
Substituting the norms $\|\bbx_N^\perp\|$  and $\|\nabla F(\bbx_N)\|$ in \eqref{eq:subopt-Cauchy} their upper bounds in \eqref{eq:xperp} and ~\eqref{upper_bound_on_primal_grad_norm} leads to
\begin{align*}
    F(\bbx_N) - F(\bbx^*) \le \Oc{\left(\sqrt{S\tfrac{\lambda_{\max}(\bbL)^3}{\mu^3\lambda^{+}_{\min}(\bbL)}}MN^{\frac{-s}{s+1}}\right)}.
\end{align*}
Hence, the claim in \eqref{main_theorem_third_claim} follows.
\end{proof}
\vspace{2mm}

Here we compare our result against known results. In terms of the consensus distance $\|\sqrt{\bbL}\bbx^*(\bby_N)\|^2$ and the distance to primal optimum $\|\bbx_N - \bbx^*\|^2$, Our proposed algorithm achieves the rate $\Oc{(N^{-2s/(s+1)})}$, which is faster than $\Oc{(N^{-1})}$ when $s \ge 2$, and approaches the optimal rate $\Oc{(N^{-2})}$ as $s \to \infty$. In our experiments presented in the following section, we observe that the rate is matched when $s=4$ and we suspect that the analysis is conservative. 

The suboptimality $F(\bbx_N) - F(\bbx^*)$ bound approaches $\Oc{(N^{-1})}$ as $s \to \infty$. We would like to emphasize two points. First, simply applying gradient descent on the dual function followed by a similar analysis will give a sublinear rate of $\Oc{(N^{-1/2})}$. Second, under additional assumptions such as $L(F(\bbx_N) - F(\bbx^*)) \ge \|\nabla F(\bbx_N)\|^2$ (which is implied by and more general than Lipschitz-gradient), \eqref{eq:subopt-Cauchy} generates a convergence rate of $\Oc{(N^{-2s/(s+1)})}$. 

Note that the distance to consensus and the distance to optimally, in terms of the optimization variable, tend to the optimal lower bounds with increasing orders of the integrator. However, explicit dependencies on the spectral properties of the graph are suboptimal. As it was shown in~\cite{uribe2018dual}, for distributed problems the optimal rates have a dependency of $\sqrt{{\lambda_{\max}(\bbL)}/{\lambda_{\min}^+(\bbL)}}$ in terms of the graph Laplacian $\bbL$. The dependency on the function parameters $\mu$ and $L$ are suboptimal as well. Achieving optimality in terms of the function and graph parameters requires further investigation. In general, one would expect that the only loss in optimality is with respect to $S$, i.e., the number of additional oracle calls, which is a constant factor. More importantly, whether the slow convergence of suboptimality is an artifact of the proof method or the discretization approach used remains an open question.

\section{Numerical experiments}\label{sec:numerical}

In this section, we present two numerical experiments to validate our theoretical results. First, we study a distributed linear regression problem, where a network of agents seeks to solve the following optimization problem:

\vspace{-0.5cm}
\begin{align}\label{ex:cost1}
    \min_{z \in \mathbb{R}^n} \frac{1}{2m}\|Hz - b\|_2^2,
    \end{align}
    where $m$ is the total number of available data points, $p$ is the dimension of the data points, $b \in \mathbb{R}^{m}$, and $H \in \mathbb{R}^{m \times p}$.
    
Following the reformulation of Problem~\eqref{ex:cost1} described in Section~\ref{sec:problem}, we can state the problem in its distributed form

\vspace{-0.5cm}
\begin{align}\label{ex:cost_dist}
    \min_{\bbx\in \reals^{np}}\sum_{i=1}^{n}  \frac{1}{2}\frac{1}{nl}\|b_i - H_ix_i\|_2^2 \qquad \text{s.t.}\ \sqrt{\bbL} \bbx=0.
\end{align}
Here, $n$ is the number of nodes, $l$ is the number of data points per node.  $b_i \in \mathbb{R}^l$ and $H_i \in \mathbb{R}^{l\times p}$ for each $i$ are the subset of points available to agent $i$. The points are generated randomly form uniform distribution.

\begin{figure}[ht]
    \centering
    \resizebox{7cm}{!}{
    \begin{tikzpicture}
    \draw (0,0) -- (7,0) -- (7,0.5) -- (0,0.5) -- (0,0);
    \draw[black,line width=3pt] (0.25,0.25) -- (0.75,0.25) node[black] at (1.2,0.25) {$\text{CGD}$};
    \draw[blue,line width=3pt] (2,0.25) -- (2.5,0.25) node[black] at (3.1,0.25) {$\text{DAGM}$};
    \draw[red,line width=3pt] (4,0.25) -- (4.5,0.25) node[black] at (4.9,0.25) {$\text{DGD}$};
    \draw[green,line width=3pt] (5.5,0.25) -- (6,0.25) node[black] at (6.5,0.25) {$\text{Alg.~\ref{alg_1}}$};
    \end{tikzpicture} } 
    \\
    \subfigure{
    \begin{tikzpicture}
        \begin{axis}[
        width=4.4cm,height=3.4cm,scale=1,
        x label style={at={(axis description cs:0.5,0.05)},anchor=north},
        y label style={at={(axis description cs:0.15,.5)},anchor=south},
        title = {$\mid F(\bbx_k) - F^* \mid$},
        ylabel={Star Graph},
        ticklabel style = {font=\tiny},
        ymode = log,
        xmode = log,
        ymin = 1e-15, ymax=1e-3,
        xmin = 4, xmax=1e4, 
        every axis plot/.append style={line width=1pt}],
        legend pos=south west;
        \addplot [blue]     table [x index=30,y index=25]{L_optimal2.dat};
        \addplot [red]         table [x index=30,y index=27]{L_optimal2.dat};
        \addplot [black]    table [x index=30,y index=24]{L_optimal2.dat};
        \addplot [green]    table [x expr={4*\thisrowno{18}},y index=16]{L_discrete_abs3.dat};
        \end{axis}
    \end{tikzpicture}
    \begin{tikzpicture}
        \begin{axis}[
        width=4.4cm,height=3.4cm,scale=1,
        x label style={at={(axis description cs:0.5,0.05)},anchor=north},
        y label style={at={(axis description cs:0.15,.5)},anchor=south},
        title = {$\|\bbL\bbx_k\|_2$},
        ticklabel style = {font=\tiny},
        ymode = log,
        xmode = log,
        ymin = 1e-8, ymax=1e1,
        xmin = 4, xmax=1e4, 
        every axis plot/.append style={line width=1pt}],
        legend pos=south west;
        \addplot [blue]     table [x index=30,y index=26]{L_optimal2.dat};
        \addplot [red]         table [x index=30,y index=28]{L_optimal2.dat};
        \addplot [green]    table [x expr={4*\thisrowno{18}},y expr={sqrt(\thisrowno{17})}]{L_discrete_abs3.dat};
        \end{axis}
    \end{tikzpicture}
    }
    \\
        \subfigure{
    \begin{tikzpicture}
        \begin{axis}[
        width=4.4cm,height=3.4cm,scale=1,
        x label style={at={(axis description cs:0.5,0.05)},anchor=north},
        y label style={at={(axis description cs:0.15,.5)},anchor=south},
        ylabel={Cycle Graph},
        ticklabel style = {font=\tiny},
        ymode = log,
        xmode = log,
        ymin = 1e-15, ymax=1e-3,
        xmin = 4, xmax=1e4, 
        every axis plot/.append style={line width=1pt}],
        legend pos=south west;
        \addplot [blue]         table [x index=30,y index=5]{L_optimal2.dat};
        \addplot [red]         table [x index=30,y index=7]{L_optimal2.dat};
        \addplot [black]             table [x index=30,y index=4]{L_optimal2.dat};
        \addplot [green]     table [x expr={4*\thisrowno{18}},y index=4]{L_discrete_abs3.dat};
        \end{axis}
    \end{tikzpicture}
    \begin{tikzpicture}
        \begin{axis}[
        width=4.4cm,height=3.4cm,scale=1,
        x label style={at={(axis description cs:0.5,0.05)},anchor=north},
        y label style={at={(axis description cs:0.15,.5)},anchor=south},
        ticklabel style = {font=\tiny},
        ymode = log,
        xmode = log,
        ymin = 1e-8, ymax=1e1,
        xmin = 4, xmax=1e4, 
        every axis plot/.append style={line width=1pt}],
        legend pos=south west;
        \addplot [blue]         table [x index=30,y index=6]{L_optimal2.dat};
        \addplot [red]         table [x index=30,y index=8]{L_optimal2.dat};
        \addplot [green]     table [x expr={4*\thisrowno{18}},y expr={sqrt(\thisrowno{5})}]{L_discrete_abs3.dat};
        \end{axis}
    \end{tikzpicture}
    }
    \\
        \subfigure{
    \begin{tikzpicture}
        \begin{axis}[
        width=4.4cm,height=3.4cm,scale=1,
        x label style={at={(axis description cs:0.5,0.05)},anchor=north},
        y label style={at={(axis description cs:0.15,.5)},anchor=south},
        ylabel={Erd\H{o}s-R\'enyi},
        ticklabel style = {font=\tiny},
        ymode = log,
        xmode = log,
        ymin = 1e-15, ymax=1e-4,
        xmin = 4, xmax=1e4, 
        every axis plot/.append style={line width=1pt}],
        legend pos=south west;
        \addplot [blue]         table [x index=30,y index=15]{L_optimal2.dat};
        \addplot [red]         table [x index=30,y index=17]{L_optimal2.dat};
        \addplot [black]             table [x index=30,y index=14]{L_optimal2.dat};
        \addplot [green]     table [x expr={4*\thisrowno{18}},y index=10]{L_discrete_abs3.dat};
        \end{axis}
    \end{tikzpicture}
    \begin{tikzpicture}
        \begin{axis}[
        width=4.4cm,height=3.4cm,scale=1,
        x label style={at={(axis description cs:0.5,0.05)},anchor=north},
        y label style={at={(axis description cs:0.15,.5)},anchor=south},
        ticklabel style = {font=\tiny},
        ymode = log,
        xmode = log,
        ymin = 1e-12, ymax=1e1,
        xmin = 4, xmax=1e4, 
        every axis plot/.append style={line width=1pt}],
        legend pos=south west;
        \addplot [blue]         table [x index=30,y index=16]{L_optimal2.dat};
        \addplot [red]         table [x index=30,y index=18]{L_optimal2.dat};
        \addplot [green]     table [x expr={4*\thisrowno{18}},y expr={sqrt(\thisrowno{11})}]{L_discrete_abs3.dat};
        \end{axis}
    \end{tikzpicture}
    }
    \caption{Distance to optimality and consensus for Problem~\eqref{ex:cost_dist} over various graphs with $100$ agents and $100$ data points for each agent.}
    \label{fig:regression}
\end{figure}
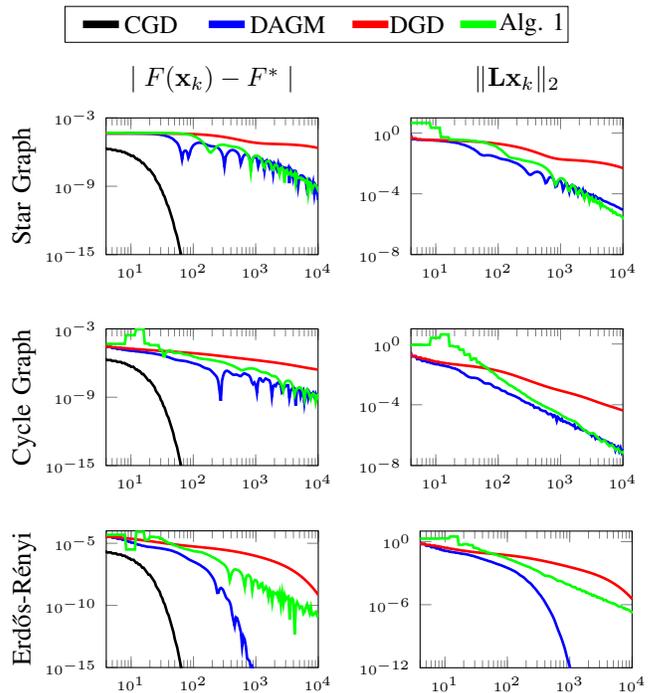 

Figure~\ref{fig:regression} shows simulation results for Problem~\eqref{ex:cost_dist} using a star graph, a cycle graph, and a Erd\H{o}s-R\'enyi random graph. Each graph has $100$ nodes and each node holds $100$ data points. Moreover, we compare the performance of our algorithm with a centralized method (CGD), the optimal distributed method (DAGM)~\cite{uribe2018dual}, and distributed gradient descent (DGD). For these simulations, we have chosen the order of the integration to be $s=4$. Results show that even with a relatively small order of the integrator, the performance of Algorithm~\ref{alg_1} is comparable with the optimal method proposed in~\cite{uribe2018dual}, both in terms of the distance to optimality and distance to consensus.

\begin{figure}[ht]
    \centering
    \resizebox{6cm}{!}{
    \begin{tikzpicture}
    \draw (0,0) -- (6,0) -- (6,0.5) -- (0,0.5) -- (0,0);
    \draw[blue,line width=3pt] (0.25,0.25) -- (0.75,0.25) node[black] at (1.3,0.25) {$\text{DAGM}$};
    \draw[red,line width=3pt] (2.5,0.25) -- (3,0.25) node[black] at (3.5,0.25) {$\text{DGD}$};
    \draw[green,line width=3pt] (4.3,0.25) -- (4.8,0.25) node[black] at (5.3,0.25) {$\text{Alg.}~1$};
    \end{tikzpicture} } 
    \\
    \subfigure{
    \begin{tikzpicture}
        \begin{axis}[
        width=4.4cm,height=3.4cm,scale=1,
        x label style={at={(axis description cs:0.5,0.05)},anchor=north},
        y label style={at={(axis description cs:0.15,.5)},anchor=south},
        title = {$\mid F(\bbx_k) - F^* \mid$},
        ylabel={Star Graph},
        ticklabel style = {font=\tiny},
        ymode = log,
        xmode = log,
        ymin = 1e-4, ymax=1e2,
        xmin = 1e1, xmax=1e4, 
        every axis plot/.append style={line width=1pt}],
        legend pos=south west;
        \addplot [blue]         table [x index=12,y index=1]{L_optimal_entropy.dat};
        \addplot [red]             table [x index=12,y index=0]{L_optimal_entropy.dat};
        \addplot [green]        table [x expr={4*\thisrowno{6}},y index=0]{L_discrete_entropy.dat};
        \end{axis}
    \end{tikzpicture}
    \begin{tikzpicture}
        \begin{axis}[
        width=4.4cm,height=3.4cm,scale=1,
        x label style={at={(axis description cs:0.5,0.05)},anchor=north},
        y label style={at={(axis description cs:0.15,.5)},anchor=south},
        title = {$\|\bbL\bbx_k\|_2$},
        ticklabel style = {font=\tiny},
        ymode = log,
        xmode = log,
        ymin = 1e-7, ymax=1e1,
        xmin = 4, xmax=1e4, 
        every axis plot/.append style={line width=1pt}],
        legend pos=south west;
        \addplot [blue]     table [x index=12,y index=3]{L_optimal_entropy.dat};
        \addplot [red]         table [x index=12,y index=2]{L_optimal_entropy.dat};
        \addplot [green]    table [x expr={4*\thisrowno{6}},y expr={sqrt(\thisrowno{1})}]{L_discrete_entropy.dat};
        \end{axis}
    \end{tikzpicture}
    }
    \\
        \subfigure{
    \begin{tikzpicture}
        \begin{axis}[
        width=4.4cm,height=3.4cm,scale=1,
        x label style={at={(axis description cs:0.5,0.05)},anchor=north},
        y label style={at={(axis description cs:0.15,.5)},anchor=south},
        ylabel={Cycle Graph},
        ticklabel style = {font=\tiny},
        ymode = log,
        xmode = log,
        ymin = 1e-7, ymax=1e2,
        xmin = 1e1, xmax=1e4, 
        every axis plot/.append style={line width=1pt}],
        legend pos=south west;
        \addplot [blue]         table [x index=12,y index=5]{L_optimal_entropy.dat};
        \addplot [red]             table [x index=12,y index=4]{L_optimal_entropy.dat};
        \addplot [green]        table [x expr={4*\thisrowno{6}},y index=2]{L_discrete_entropy.dat};
        \end{axis}
    \end{tikzpicture}
    \begin{tikzpicture}
        \begin{axis}[
        width=4.4cm,height=3.4cm,scale=1,
        x label style={at={(axis description cs:0.5,0.05)},anchor=north},
        y label style={at={(axis description cs:0.15,.5)},anchor=south},
        ticklabel style = {font=\tiny},
        ymode = log,
        xmode = log,
        ymin = 1e-6, ymax=1e1,
        xmin = 4, xmax=1e4, 
        every axis plot/.append style={line width=1pt}],
        legend pos=south west;
        \addplot [blue]         table [x index=12,y index=7]{L_optimal_entropy.dat};
        \addplot [red]             table [x index=12,y index=6]{L_optimal_entropy.dat};
        \addplot [green]        table [x expr={4*\thisrowno{6}},y expr={sqrt(\thisrowno{3})}]{L_discrete_entropy.dat};
        \end{axis}
    \end{tikzpicture}
    }
    \\
        \subfigure{
    \begin{tikzpicture}
        \begin{axis}[
        width=4.4cm,height=3.4cm,scale=1,
        x label style={at={(axis description cs:0.5,0.05)},anchor=north},
        y label style={at={(axis description cs:0.15,.5)},anchor=south},
        ylabel={Erd\H{o}s-R\'enyi},
        ticklabel style = {font=\tiny},
        ymode = log,
        xmode = log,
        ymin = 1e-4, ymax=1e2,
        xmin = 1e2, xmax=1e4, 
        every axis plot/.append style={line width=1pt}],
        legend pos=south west;
        \addplot [blue]         table [x index=12,y index=9]{L_optimal_entropy.dat};
        \addplot [red]             table [x index=12,y index=8]{L_optimal_entropy.dat};
        \addplot [green]        table [x expr={4*\thisrowno{6}},y index=4]{L_discrete_entropy.dat};
        \end{axis}
    \end{tikzpicture}
    \begin{tikzpicture}
        \begin{axis}[
        width=4.4cm,height=3.4cm,scale=1,
        x label style={at={(axis description cs:0.5,0.05)},anchor=north},
        y label style={at={(axis description cs:0.15,.5)},anchor=south},
        ticklabel style = {font=\tiny},
        ymode = log,
        xmode = log,
        ymin = 1e-4, ymax=1e1,
        xmin = 4, xmax=1e4, 
        every axis plot/.append style={line width=1pt}],
        legend pos=south west;
        \addplot [blue]         table [x index=12,y index=11]{L_optimal_entropy.dat};
        \addplot [red]             table [x index=12,y index=10]{L_optimal_entropy.dat};
        \addplot [green]        table [x expr={4*\thisrowno{6}},y expr={sqrt(\thisrowno{5})}]{L_discrete_entropy.dat};
        \end{axis}
    \end{tikzpicture}
    }
    \caption{Distance to optimality and consensus for Problem~\eqref{problem:entropy} over various graphs with $100$ agents with a distribution over $100$ points for each agent.}
    \label{fig:entropy}
\end{figure}
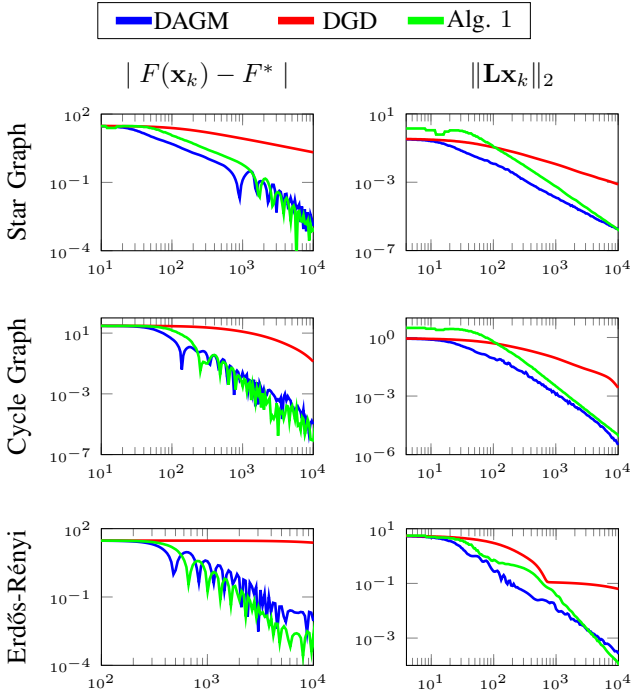

As a second example, we consider the Kullback-Leibler~(KL) barycenter computation problem \cite{Uribe2018}. This problem is strongly convex and $M$-Lipschitz, which is defined as
    \begin{align}\label{problem:entropy}
    \min_{z \in S_p(1)} \sum\limits_{i=1}^{n}D_{KL}(z\|q_i) \triangleq \sum\limits_{i=1}^{n}\sum\limits_{j=1}^{p}z^i \log\left({z_i}/{[q_i]_j} \right),  
    \end{align}
where \mbox{$S_p(1) = \{ z \in \mathbb{R}^p: z_j \geq 0 ; j=1,\dots,n; \sum_{j=1}^{p}z_j = 1 \}$} is a unit simplex in $\mathbb{R}^p$ and $q_i \in S_p(1)$ for all $i$. Each agent has a private probability distribution $q^i$ and seeks to compute the a probability distribution that minimizes the average KL distance to the distributions $\{q_i\}_{i=1}^n$. Figure \ref{fig:entropy} shows the results for the KL barycenter problem for a cycle, a star and an Erd\H{o}s-R\'enyi random graph with $n=100$, and $p=100$. We show the distance to optimality as well as the distance to consensus.

In our first experiment, DAGM and DGD eventually achieve a linear rate because the $L_2$ norm is smooth and strongly convex. Algorithm~\ref{alg_1}, however, achieves a polynomial rate. In the second example, since the KL problem is not smooth, all algorithms converge at polynomial rates. In this case, Algorithm~\ref{alg_1} achieves the best relative performance.

\begin{figure}[tp]
    \centering
    \subfigure{
    \begin{tikzpicture}
        \node at (0.7,1.5) {{\footnotesize {\color{blue} \textbf{--}} $s=1$}};
		\node at (0.7,1.2) {{\footnotesize {\color{red} \textbf{--}} $s=2$}};
		\node at (0.7,0.9) {{\footnotesize {\color{green} \textbf{--}} $s=4$}};
        \begin{axis}[
        width=4.4cm,height=3.4cm,scale=1,
        x label style={at={(axis description cs:0.5,0.05)},anchor=north},
        y label style={at={(axis description cs:0.15,.5)},anchor=south},
        title = {$\mid F(\bbx_k) - F^* \mid$},
        ylabel={Erd\H{o}s-R\'enyi},
        ticklabel style = {font=\tiny},
        ymode = log,
        xmode = log,
        ymin = 1e-7, ymax=1e2,
        xmin = 4, xmax=1e4, 
        every axis plot/.append style={line width=1pt}],
        legend pos=south west;
        \addplot [blue]         table [x expr={4*\thisrowno{6}},y index=0]{L_discrete_entropy_s.dat};
        \addplot [red]             table [x expr={4*\thisrowno{6}},y index=2]{L_discrete_entropy_s.dat};
        \addplot [green]        table [x expr={4*\thisrowno{6}},y index=4]{L_discrete_entropy_s.dat};
        \end{axis}
    \end{tikzpicture}
    \begin{tikzpicture}
        \begin{axis}[
        width=4.4cm,height=3.4cm,scale=1,
        x label style={at={(axis description cs:0.5,0.05)},anchor=north},
        y label style={at={(axis description cs:0.15,.5)},anchor=south},
        title = {$\|\bbL\bbx_k\|_2$},
        ticklabel style = {font=\tiny},
        ymode = log,
        xmode = log,
        ymin = 1e-7, ymax=1e1,
        xmin = 1e2, xmax=1e4, 
        every axis plot/.append style={line width=1pt}],
        legend pos=south west;
        \addplot [blue]         table [x expr={4*\thisrowno{6}},y index=1]{L_discrete_entropy_s.dat};
        \addplot [red]             table [x expr={4*\thisrowno{6}},y index=3]{L_discrete_entropy_s.dat};
        \addplot [green]        table [x expr={4*\thisrowno{6}},y index=5]{L_discrete_entropy_s.dat};
        \end{axis}
    \end{tikzpicture}
    }
    \caption{The effect of increasing the order of the integrator.}
    \label{fig:s}
\end{figure}
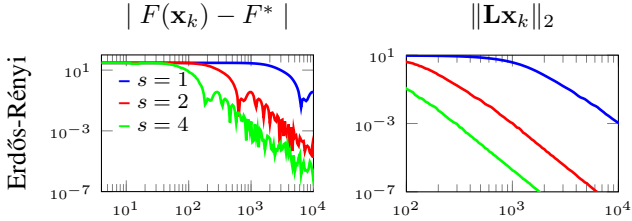

Finally, Figure~\ref{fig:s} shows how the order of the integrator affects the convergence rate of Algorithm~\ref{alg_1}. We test $s=1$, $s=2$, and $s=4$, and as $s$ increases we observe faster rates.

\section{Conclusions and Future Work}\label{sec:conclusions}

We proposed a new distributed accelerated algorithm that achieves a faster convergence rate than gradient descent. The algorithm follows a simple intuition by directly discretizing the heavy-ball ODE on the appropriately formulated dual problem. This approach demonstrates that tools and results from dynamical system theory can be applied to optimization and provide insights into existing problems.

The proposed method requires an exact solution of the inner maximization problem. One can study the effects of having approximate solutions only, but this is left for future work. Additionally, we point out that the convergence rate estimates resulting from our analysis are strictly sub-optimal. Analyze convergence under other weaker convexity or smoothness assumptions also requires further study.

\bibliographystyle{IEEEtran}
\bibliography{references,all_refs}
\end{document}